\documentclass{article}

\title{Cyclic coverings and Seshadri constants on smooth surfaces.}
\author{Luis Fuentes Garc\'{\i}a}
\date{}
\newtheorem{teo}{Theorem}[section]

\newtheorem{prop}[teo]{Proposition}
\newtheorem{cor}[teo]{Corollary}
\newtheorem{lemma}[teo]{Lemma}

\newtheorem{rem}[teo]{Remark}
\newtheorem{conj}[teo]{Conjecture}
\newtheorem{question}[teo]{Question}

\def\mod{\mathop{\rm mod}}

\def\mod{\mathop{\rm mod}}

\def\g2{\pi}

\def\P{{\bf P}}

\newcommand\Te{{\cal O}}

\def\ZZ{\leavevmode\hbox{$\rm Z$}}

\def\impp{{\quad\Rightarrow\quad}}

\def\qed{\hspace{\fill}$\rule{2mm}{2mm}$}

\newcommand\lrw{\longrightarrow}

\begin{document}

\maketitle

\vspace{0.1cm}

\begin{abstract}

We study the Seshadri constants of cyclic coverings of smooth surfaces. The existence of an automorphism on these surfaces can be used to produce Seshadri exceptional curves. We give a bound for multiple Seshadri constants on cyclic coverings of surfaces with Picard number $1$. Morevoer, we apply this method to $n$-cyclic coverings of the projective plane. When $2\leq n\leq 9$, explicit values are obtained. We relate this problem with the Nagata conjecture. 

{\bf MSC (2000):} Primary 14E20; secondary, 14C20.

 {\bf Key Words:} Cyclic coverings, Seshadri constants, Nagata conjecture.

\end{abstract}

\section{Introduction.}

The multiple Seshadri constants are a natural generalization of the Seshadri constants at single points defined by Demailly in \cite{De92}. If $X$ is a smooth projective surface, $L$ is a nef bundle on $X$ and $P_1,\ldots,P_r$ are distinct points in $X$, then the Seshadri constant of $L$ at $P_1,\ldots,P_r$ is::
$$
\epsilon(L; P_1,\ldots,P_r)=inf \frac{C\cdot L}{\sum_{i=1}^r mult_{P_i}C},
$$
where $C$ runs over all curves passing through at least one of the points $P_1,\ldots,P_r$. These constants have the upper bound:
$$
\epsilon(L; P_1,\ldots, P_r)\leq \sqrt{\frac{L^2}{r}}.
$$
When the constant does not reach this bound, there is a curve with high multiplicity at the points such that 
$$
\epsilon(L; P_1,\ldots,P_r )=\frac{C\cdot L}{\sum_{i=1}^r mult_{P_i}C}.
$$
These curves are called Seshadri exceptional curves (see \cite{EiKuLa95}). 

When $r=1$, general bounds for the Seshadri constants on surfaces are given in \cite{Ba99}, \cite{Na03} or \cite{St98}. An interesting open problem is their irrationality. It seems that surfaces with irrational Seshadri constants must exist. However, the  explicit known values are always rational. They were computed for simple abelian surfaces by Th. Bauer and T. Szemberg (see \cite{Ba98}); Ch. Schultz gave values for Seshadri constants on products of two elliptic curves (see \cite{Sc04}); I obtained the Seshadri constants on elliptic ruled surfaces (see \cite{Fu06}). Note that, in all these cases, the Seshadri exceptional curves providing the constants are built in the same way. The existence of an involution on these surfaces allows to construct invariant curves with high multiplicity at the fixed points.

In this paper we use this idea to study Seshadri constants on $n$-cyclic coverings of smooth surfaces:
$$
\pi:X\lrw Y.
$$
These surfaces have an automorphism of order $n$. We search Seshadri exceptional curves on invariant linear systems. With this method, we will obtain results about the Seshadri constants of line bundles $\pi^*L$ of $X$ at points on the ramification divisor.

We extend a result of Steffens about simple Seshadri constants on surfaces with Picard number $1$ (see \cite{St98}), for multiple Seshadri constants on cyclic coverings of smooth surfaces. We prove :

\begin{teo}

Let $\pi:X\lrw Y$  be a branched $n$-cyclic covering of a smooth surface $Y$ with $\rho(Y)=1$. Let $L$ be an ample generator of $NS(Y)$. Then, if $P_1,\ldots,P_r$ are very general points on $X$ then:
$$
\frac{\left[\sqrt{r(\pi^*L)^2}\right]}{r}\leq \epsilon(\pi^*L;P_1,\ldots,P_r)\leq \frac{\sqrt{(\pi^*L)^2}}{\sqrt{r}}.
$$
In particular, if $\sqrt{rnL^2}$ is an integer, then $\epsilon(\pi^*L;P_1,\ldots,P_r)=\frac{\sqrt{(\pi^*L)^2}}{\sqrt{r}}$.

\end{teo}

Some results in the same direction have been obtained for B. Harbourne in \cite{Ha03}. However, the he needs $L$ to be very ample and $r\geq \pi^*L^2$. In \cite{Ro041}, J. Roe proves a similar bound, but he uses the Nagata conjecture. 

We also apply the method to study $n$-cyclic coverings of the projective plane. When $2\leq n\leq 9$ we obtain explicit values of the Seshadri constant of $\pi^*\Te_{P^2}(1)$ at a very general point $x$ on the ramification divisor (Theorem \ref{nplane}):
$$
\begin{array}{|c|c|c|c|c|c|c|c|c|}
\hline
{n}&{2}&{3}&{4}&{5}&{6}&{7}&{8}&{9}\\
\hline
{\epsilon(\pi^*\Te_{P^2}(1),x)}&{1}&{3/2}&{2}&{2}&{12/5}&{21/8}&{48/17}&{3}\\
\hline
\end{array}
$$

Finally, we see the relation between the study of the Seshadri constant of $\pi^*\Te_{P^2}(1)$ and the Nagata conjecture. In particular, the Seshadri exceptional curves for $\pi^*\Te_{P^2}(1)$ are given by the pullback of curves in $\P^2$ passing through $n$ infinitely near points with prescribed order. In fact, this is true for any covering (see Proposition \ref{clusters}). In general, we will apply this relation to prove:

\begin{teo}

Let $\pi:X\lrw Y$  be a branched $n$-cyclic covering of a smooth surface $Y$. Let $L$ be an ample divisor of $NS(Y)$. Then, if $P_1,\ldots,P_r$ are different points in the branch divisor of $\pi$:
$$
\epsilon(\pi^*L;P_1,\ldots,P_r)\leq n\epsilon(L; n\, r)
$$
where $\epsilon(L; n\, r)$ is the Seshadri constant of $L$ at $n\cdot r$ very general points.

\end{teo}

We refer to \cite{Ba99} for a systematic study of the main properties of the Seshadri constants on surfaces. 

{\bf Acknowledgments: } I thank A. Broustet for his remarks about Seshadri constants on Del Pezzo surfaces.

\section{Cyclic coverings.}

Firstly, we recall some well known facts about cyclic coverings (see \cite{BaPeVa84}). Let $Y$ be a smooth surface and let $\Te_X(M)$ be a line bundle on $Y$.  Consider a smooth reduced divisor $B\sim nM$ on $Y$. Then we have a $n$-cyclic covering:
$$
\pi:X\lrw Y
$$
with branch divisor $B$. Since $B$ is smooth and reduced, $X$ is a smooth surface. Let $R$ be the reduced divisor $\pi^{-1}(B)$ (the ramification divisor). It holds:

\begin{enumerate}

\item $\Te_X(R)=\pi^*\Te_Y(M)$.

\item $\pi^*B=nR$.

\item $\Te_X(K_X)=\pi^*(\Te_Y(K_Y+(n-1)M))$.

\end{enumerate}

On the other hand, there is an induced automorphism of order $n$:
$$
\sigma:X\lrw X.
$$
The fixed points of $\sigma$ are exactly the points of the ramification divisor $R$. Moreover, we have automorphisms:
$$
\bar \sigma:H^0(X,\pi^*L)\lrw H^0(X,\pi^*L), \quad \hbox { where $L$ is a line bundle on $Y$.}
$$
The spaces of eigenvectors of these automorphisms are given by the following decomposition:
$$
H^0(X,\pi^*L)\cong \bigoplus_{k=0}^{n} H^0(Y,L-kM).
$$
In fact, if  $D\sim \pi^*L$ is a $\sigma$-invariant divisor then $D=E+kR$, where $E\sim \pi^*(L-kM)$. Since we are interested on irreducible divisors, we will study invariant divisors associated to the eigenvalue $1$.

We will denote by $H^0(X,\pi^*L)_1$ the space of sections associated to the eigenvalue $1$. Let us consider a point $x$ at the ramification divisor $R$. We study the existence of divisors in $H^0(X,\pi^*L)_1$ passing through $x$ with given multiplicity. 

Let us take local coordinates $(u,v)$ such that $v=0$ is the local equation of $R$ and $u=0$ corresponds to an invariant divisor passing through $x$. The automorphism $\sigma$ have the following local expression:
$$
\sigma(u,v)=(u,\theta v),
$$
where $\theta$ is a primitive $n$-root of unity.

Let $f(u,v)=0$ be a local equation of a divisor in $H^0(X,\pi^*L)_1$. It must verify:
$$
f(u,v)=f(u,\theta v).
$$
From this:
$$
\frac{\partial f}{\partial u^i \partial v^j}(0,0)=\theta^j\frac{\partial f}{\partial u^i \partial v^j}(0,0)
$$
and then:
$$
\frac{\partial f}{\partial u^i \partial v^j}(0,0)=0,\quad \hbox{ when }j\neq 0 \mod n.
$$
Therefore, we deduce the following lemma:

\begin{lemma} \label{multiplicity}

With the previous notation, the number of conditions on a divisor in $H^0(X,\pi^*L)_1$ to pass through $x\in R$ with multiplicity at least $m$ is:
$$
(k+1)\left({\frac{nk}{2}+r}\right),\quad \hbox{ where } m=nk+r, \quad 0\leq r\leq n-1.
$$

\end{lemma} \qed

\section{Cyclic coverings of smooth surfaces with Picard number $1$.}

In \cite{St98}, Steffens proved the following result:

\begin{teo}\label{Steffens}

Let $X$ be a surface with $\rho(X)=rk(NS(X))=1$ and let $L$ be an ample generator of $NS(X)$. Let $\alpha$ be an integer with $\alpha^2\leq L^2$. If $\eta\in X$ is a very general point, then $\epsilon(L,\eta)\geq \alpha$. In particular, if $\sqrt{L^2}$ is an integer, then $\epsilon(L,\eta)=\sqrt{L^2}$.

\end{teo}

The proof is based on two facts. First, the following result of Ein-Lazarsfeld (\cite{EiLa93}):

\begin{lemma}\label{Ein}

Let $\{C_t\in x_t\}_{t\in\Delta}$ be a $1$-parameter family of reduced irreducible curves on a smooth projective surface $X$, such that $mult_{x_t}(C_t)\geq m$ for all $t\in \Delta$. Then:
$$
(C_t)^2\geq m(m-1).
$$
\end{lemma}

Moreover, because $\rho(X)=1$, he can use that a Seshadri exceptional curve $C$ is numerically equivalent to $dL$ for some integer $d$.

This result can be generalized to multiple Seshadri constants of cyclic coverings of surfaces with Picard number $1$. Note that a cyclic covering
$
\pi:X\lrw Y
$
verifies $\rho(X)\geq \rho(Y)$.

We will need two previous lemmas which are useful to extend the Lemma \ref{Ein}:

\begin{lemma}\label{Xu}

Let $X$ be a smooth projective surface, $\{C_t,y_t\}_{t\in \Delta}$ a non-trivial one parameter
family consisting of integral $C_t\subset X$ and $y_1,\ldots,y_r\in X$ be distinct points different from $y_t$ for at least one $t$ such that $mult_{y_i}(C_t)\geq m_i$ with $mult_{y_i}(C_t)\geq m_i$ for al $t\in \Delta$ and $i=1,\ldots,s$. If moreover $mult_{y_t}(C_t)\geq m> 0$ for all $t\in \Delta$, then:
$$
C_t^2\geq m(m-1)+\sum_{i=1}^r m_i^2.
$$

\end{lemma}
{\bf Proof:} See \cite{Ku96}. \qed

\begin{lemma}\label{numerico}

Let $m_1,\ldots,m_r$ be non negative integers. Let $M=\displaystyle\sum_{i=1}^r m_i$. Then:
$$
r(\sum_{i=1}^r m_i^2-m_{i_0})\geq M(M-1)
$$
for any integer $i_0$, with $1\leq i_0\leq r$.

\end{lemma}
{\bf Proof:} 
$$
\begin{array}{rl}
{r\displaystyle\sum_{i=1}^r m_i^2-rm_{i_0}}&{=(\displaystyle\sum_{i=1}^r m_i)^2+\displaystyle\sum_{1\leq i<j\leq r} (m_i-m_j)^2-rm_{i_0}=}\\
{}&{=M(M-1)+\displaystyle\sum_{1\leq i<j\leq r} (m_i-m_j)^2+\displaystyle\sum_{i=1}^r (m_i-m_{i_0})\geq}\\
{}&{\geq M(M-1)+\displaystyle\sum_{i\neq i_0} (m_i-m_{i_0})^2+\displaystyle\sum_{i\neq i_0}(m_i-m_{i_0})\geq}\\
{}&{\geq M(M-1).}\\
\end{array} 
$$ \qed

Now, we can prove the extended version of the Theorem \ref{Steffens}.

\begin{teo}\label{picard1}

Let $\pi:X\lrw Y$  be a branched $n$-cyclic covering of a smooth surface $Y$ with $\rho(Y)=1$. Let $L$ be an ample generator of $NS(Y)$. Then, if $P_1,\ldots,P_r$ are very general points on $X$ then:
$$
\frac{\left[\sqrt{r(\pi^*L)^2}\right]}{r}\leq \epsilon(\pi^*L;P_1,\ldots,P_r)\leq \frac{\sqrt{(\pi^*L)^2}}{\sqrt{r}}.
$$
In particular, if $\sqrt{rnL^2}$ is an integer, then $\epsilon(\pi^*L;P_1,\ldots,P_r)=\frac{\sqrt{(\pi^*L)^2}}{\sqrt{r}}$.

\end{teo}
{\bf Proof:} Since the Seshadri constants are lower semi-continuous (see \cite{Og02}), it is sufficient to prove the result for a particular set of $r$ points. Let $\sigma$ be the automorphism of $X$ induced by $\pi$. Ler $R$ be the ramification divisor. Let $$Z=\{P_{1,0},\ldots,P_{1,n-1},\ldots,P_{k,0},\ldots,P_{k,n-1},Q_1,\ldots,Q_s\}$$ be a set of $r=nk+s$ points invariant by $\sigma$ with $1\leq s\leq n$ and such that:
$$
P_{i,j}\not\in R, \quad \sigma(P_{i,j})=P_{i,j+1}, \quad \hbox{with  }j\in \ZZ_{n}, \hbox{ and } Q_1,\ldots,Q_s\in R.
$$
Let $C$ be a curve providing the Seshadri constant of $\pi^*L$ at these points. Let $l$ be the smaller integer verifying, $\sigma^l(C)=C$. Then the vector of multiplicities of $C$ at the points $P_{i,0},\ldots,P_{i,n-1}$ is:
$$
(m_{i,1},\ldots,m_{i,l},m_{i,1},\ldots,m_{i,l},\ldots,m_{i,1},\ldots,m_{i,l})
$$
Let $(n_1,\ldots,n_s)$ be the multiplicities of $C$ at $Q_1,\ldots,Q_s$. We can move the point $Q_1$ along the ramification divisor $R$. We have a one-parameter family of Seshadri exceptional curves $\{C_t,(Q_1)_t\}_{t\in R}$. Note that the ramification curve is not Seshadri exceptional for $\pi^*L$ at the points of $Z$, because:
$$
\frac{R\cdot \pi^*L}{\sum_{P\in Z}mult_P(R)}= \frac{R\cdot \pi^*L}{s}\geq \frac{\pi^*L^2}{s}\geq \frac{\sqrt{\pi^*L^2}}{r}.
$$
Thus the family of curves is not trivial and we can apply the Lemma \ref{Xu} to obtain the following bound for the generic curve $C$ of the family:
\begin{equation}\label{des0}
C^2\geq \sum_{i=1}^k\sum_{j=1}^l \frac{n}{l}m_{i,j}^2+\sum_{i=1}^s n_i^2-n_1
\end{equation}
Now, consider the divisor:
$$
D=C+\sigma(C)+\sigma^2(C)+\ldots+\sigma^{l-1}(C)
$$
It is clear that:
$$
mult_{P_{i,j}}{D}=\sum_{\mu=1}^l m_{i,\mu},\quad mult_{Q_i}(D)=l \, mult_{Q_i}(C),
$$
and
$$
\sum_{P\in Z} mult_P(D)=l\sum_{P\in Z} mult_P(C).
$$
Moreover,
$$
D\cdot \pi^*L=lC\cdot \pi^*L.
$$
From this, $D$ is a divisor providing the Seshadri constant of $\pi^*L$ at the points of $Z$. It verifies:
$$
\begin{array}{rl}
{D^2}&{=lC^2+2\displaystyle\sum_{0\leq \mu<\nu< l}\sigma^{\mu}(C)\cdot \sigma^{\nu}(C)\geq }\\
{}&{}\\
{}&{\geq lC^2+2\displaystyle\sum_{i=1}^k\displaystyle\sum_{0\leq \mu<\nu< l} n m_{i,\mu}m_{i,\nu}+l(l-1)\displaystyle\sum_{i=1}^s n_i^2\geq}\\
{}&{}\\
{}&{\stackrel{(\ref{des0})}{\geq} n\displaystyle\sum_{i=1}^k (m_{i,1}+\ldots+m_{i,l})^2+\displaystyle\sum_{i=1}^s (ln_i)^2-l n_1\geq}\\
{}&{}\\
{}&{\geq \sum_{P\in Z}mult_P(D)^2-mult_{Q_1}(D)}\\
\end{array}
$$
Applying the Lemma \ref{numerico} we deduce that:
\begin{equation}\label{des1}
rD^2 \geq M(M-1) 
\end{equation}
where $M=\sum_{P\in Z}mult_P(D)$.

Now we can rescue the arguments of Steffens. Let $\alpha$ be an integer with $\alpha^2\leq r\pi^*L^2$. Suppose that:
$$
\epsilon(\pi^*L; Z)<\frac{\alpha}{r}\leq \frac{\sqrt{\pi^*L^2}}{\sqrt{r}}.
$$
Thus:
\begin{equation}\label{des2}
\frac{D\cdot \pi^*L}{M}<\frac{\alpha}{r}\impp rD\cdot \pi^*L<\alpha M.
\end{equation}
On the other hand, since $D$ is invariant and $\rho(Y)=1$ there is an integer $d$ such that $D\equiv d\pi^*L$ and:
\begin{equation}\label{des3}
\alpha^2\leq r\pi^*L^2\impp \alpha^2d\leq rd\pi^*L^2=rD\cdot \pi^*L^2\stackrel{(\ref{des2})}{<}\alpha\, M\impp \alpha d\leq M-1.
\end{equation}
Finally:
$$
M(M-1)\stackrel{(\ref{des1})}{\leq} rD^2=rdD\cdot \pi^*L\stackrel{(\ref{des2})}{<}\alpha dM\stackrel{(\ref{des3})}{\leq} M(M-1)
$$
and this is a contradiction. \qed

\begin{rem}

The argument used in this proof give us an interesting consequence. In order to compute the Seshadri constant of line bundles $\pi^*L$ at a set of invariant points of the induced automorphism $\sigma$, we only have to consider divisors that are invariant for $\sigma$. Note, that this does not depend on the Picard number of $Y$.

\end{rem}

\begin{cor}

Let $X$  be a smooth surface $X$ with $\rho(X)=1$. Let $L$ be an ample generator of $NS(X)$. If $P_1,\ldots,P_r$ are very general points on $X$ then:
$$
\frac{\left[\sqrt{rL^2}\right]}{r}\leq \epsilon(L;P_1,\ldots,P_r)\leq \frac{\sqrt{L^2}}{\sqrt{r}}.
$$
In particular, if $\sqrt{rL^2}$ is an integer, then $\epsilon(L;P_1,\ldots,P_r)=\frac{\sqrt{L^2}}{\sqrt{r}}$. \qed
\end{cor}

\section{Cyclic coverings of the projective plane.}

We will work with a $n$-cyclic covering of $\P^2$:
$$
\pi:X\lrw \P^2
$$
with branch divisor $B\sim nb L$ and $L=\Te_{P^2}(1)$. We are interested on the Seshadri constant of $\pi^*L$. A  direct application of the Theorem \ref{picard1} gives:

\begin{teo}\label{aproximacion}

Let $\pi:X\lrw \P^2$  be an $n$-cyclic covering of $\P^2$ and $L=\Te_{P^2}(1)$. If $\eta$ is a very general point on $X$ then:
$$
\left[\sqrt{n}\right]\leq \epsilon(\pi^*L,\eta)\leq \sqrt{n}.
$$
In particular, if $\sqrt{n}$ is an integer, $\epsilon(\pi^*L,\eta)=\sqrt{n}$. \qed

\end{teo} 

We can obtain a refinement of this result. Let us consider the linear system $|d\pi^*L|$. Let $x$ be a point on the ramification divisor $R$. We will use the Lemma \ref{multiplicity} to find divisors in $|d\pi^*L|_1$ with high multiplicity $m$ at $x$.

A Seshadri exceptional curve $D\sim d\pi^*L$ passing through $x$ must verify:

\begin{equation}\label{desigualdad1}
D^2\leq m^2, \quad \hbox { or equivalently, } \quad d^2n\leq m^2. 
\end{equation}

The dimension of $H^0(X,d\pi^*L)_1$ is:
$$
h^0(\Te_{P^2}(d))={{d+2}\choose 2}.
$$
This dimension must be greater than the number of conditions on a divisor $D$ to pass through $x$ with multiplicity $m$. Applying the Lemma \ref{multiplicity}, this means:
\begin{equation}\label{desigualdad2}
{{d+2}\choose 2}>(k+1)(\frac{nk}{2}+r)\iff {{d+2}\choose 2}-1\geq\left(\frac{m-r}{n}+1\right)\left(\frac{m+r}{2}\right)
\end{equation}
where $m=nk+r$, with $0\leq r<n$. From this:
$$
d^2+3d\geq\frac{m^2-r^2}{n}+m+r\stackrel{(\ref{desigualdad1})}{\impp} 3d\geq m+r-\frac{r^2}{n}\geq m.
$$
Taking squares and applying the inequality (\ref{desigualdad1}):
$$
9d^2\geq m^2\geq nd^2\impp n\leq 9.
$$
Thus, if $n>9$ we can not get the desired divisor $D$. On the other hand, if $2\leq n\leq 9$ we can expect to find values of $m$ and $d$ satisfying the inequalities (\ref{desigualdad1}) and (\ref{desigualdad2}). Explicitly, we obtain:

\begin{table}[h]
\centering
\begin{tabular}{ccccc}
{n}&{d}&{m}&{$h^0(\Te_{P^2}(d))$}&{\hbox{conditions}}\\
\hline
{2}&{1}&{2}&{3}&{2}\\
\hline
{3}&{1}&{2}&{3}&{2}\\
\hline
{4}&{1}&{2}&{3}&{2}\\
\hline
{5}&{2}&{5}&{6}&{5}\\
\hline
{6}&{2}&{5}&{6}&{5}\\
\hline
{7}&{3}&{8}&{10}&{9}\\
\hline
{8}&{6}&{17}&{28}&{27}\\
\hline
{9}&{3}&{9}&{10}&{9}\\
\hline
\end{tabular}
\caption{Seshadri exceptional divisors for $\pi^*\Te_{P^2}(1)$.}\label{tabla1}
\end{table}

Now, let us see that these exceptional (possibly reduced) divisors are optimal. We will prove two lemmas for bounding the multiplicity of any Seshadri exceptional divisor for $\pi^*L$. 

\begin{lemma}\label{cota}
Let $x$ be a generic point on the ramification divisor. Let $C\sim j\pi^*L$ be a Seshadri exceptional divisor for $\pi^*L$ at $x$. Then $mult_x(C)<h^0(\Te_{P^2}(j))$.
\end{lemma}
{\bf Proof:} Let $m=mult_x(C)$. Note that:
$$
(dL\cdot B)_{\pi(x)}=(C\cdot R)_x\geq m.
$$
This means that there is a plane curve $C'$ of degree $j$ meeting $B$ with multiplicity at least $m$ at $\pi(x)$. Let us consider the Veronesse map of degree d:
$$
v_j:\P^2\lrw \P^{h^0(\Te_{P^2}(j))-1}.
$$
Now, the curve $C'$ corresponds to a hyperplane meeting $v_j(B)$ at $y=v_j(\pi(x))$ with multiplicity at least $m$. When $y$ is a generic point, this multiplicity is upper bounded by the dimension of the ambient space. \qed

\begin{lemma}
When $n\neq 8$, the exceptional divisors related on the Table \ref{tabla1} provide the Seshadri constant of $\pi^*L$ at a very general point on the ramification divisor.
\end{lemma}
{\bf Proof:} If $n=4,9$ the result follows from the Theorem \ref{aproximacion}. 

Let $D_n$ be the Seshadri exceptional divisor described on the table \ref{tabla1}. Let $C$ be a Seshadri exceptional curve providing the Seshadri constant of $\pi^*L$ at a very general point $x$ on the ramification divisor $R$. With the same argument of the proof of the Theorem \ref{picard1}, we can construct a divisor $D\sim j\pi^*L$ invariant by the involution and providing the Seshadri constant of $\pi^*L$. The number $j$ is the smaller integer such that $|j\pi^*L|$ contains the curve $C$. Thus, since $D_n$ and $D$ are Seshadri exceptional divisors multiple of $\pi^*L$, it holds $j\leq d$. Moreover, if $m=mult_x(D)$:
\begin{equation}\label{desigualdad3}
m^2\geq D^2\geq m(m-1), \quad \hbox{ where }D^2=nj^2.
\end{equation}

Now, we use this inequality and the Lemma \ref{cota} to discard the cases which do not appear on the Table \ref{tabla1}.

\begin{enumerate}

\item If $n=2,3$, then $j\leq 1$. The unique possibility is $m=2$.

\item If $n=5,6$, then $j\leq 2$. If $j=1$, by the Lemma \ref{cota}, $m<3$ and the inequality (\ref{desigualdad3}) fails. If $j=2$, necessarily $m=5$.

\item If $n=7$, then $j\leq 3$. If $j=1$ or $j=2$, then $m<3$ or $m<6$ respectively. In both cases the inequality (\ref{desigualdad3}) fails. If $j=3$, necessarily  $m=8$. \qed

\end{enumerate}

\begin{rem}

Let us try to apply the same arguments when $n=8$. In this case $j\leq 6$. With the inequality (\ref{desigualdad3}) and Lemma \ref{cota}, we can discard the cases $j=1,2,4,5$. The problem appears when $j=3$ and $m=9$. Let us see that we can eliminate this possibility if we work with a generic branch divisor.

\end{rem}

\begin{lemma}

When $n=8$ and the branch divisor $B$ is a generic plane curve of degree $8b$, the Seshadri constant of $\pi^*L$ at a very general point $x$ on the ramification divisor is given by the divisor described on the table \ref{tabla1}.

\end{lemma}
{\bf Proof:} We must discard the case $j=3$ and $m=9$. The existence of a divisor $D\sim 3\pi^*L$ passing through $x$
 with multiplicity $9$ implies the existence of a cubic curve with multiplicity $2$ at $x$ and meeting the branch divisor with multiplicity at least $9$. Let us see that this curve does not exist. Since we suppose that $B$ is generic, it is sufficient to find an example.
 
Let us consider affine coordinates $(x,y)$. Let us take the curve $B$ given by the equation:
$$
y=x^{8b}+x^4+x^2.
$$
Now, with a direct computation, we can check that there are not any curve of degree $3$ with a singular point at $(0,0)$ and meeting $B$ with multiplicity $9$ at the same point. \qed

As a consequence of the previous discussion we have proved the following Theorem:

\begin{teo}\label{nplane}

Let $\pi:X\lrw \P^2$ be a $n$-cyclic covering of $\P^2$, with $2\leq n\leq 9$ and a generic branch divisor of degree $bn$. Let $L\sim \Te_{P^2}(1)$. Then:

\begin{enumerate}

\item If $x$ is a very general point on the ramification divisor:
$$
\begin{array}{|c|c|c|c|c|c|c|c|c|}
\hline
{n}&{2}&{3}&{4}&{5}&{6}&{7}&{8}&{9}\\
\hline
{\epsilon(\pi^*L,x)}&{1}&{3/2}&{2}&{2}&{12/5}&{21/8}&{48/17}&{3}\\
\hline
\end{array}
$$

\item If $\eta$ is a very general point on $X$:
$$
\epsilon(\pi^*L,x)\leq \epsilon(\pi^*L,\eta)\leq \sqrt{n}
$$ 

\end{enumerate}

In particular, the global Seshadri constant $\epsilon(\pi^*L)$ is rational. \qed
\end{teo}

\begin{rem}

When $n=2$ the surface $X$ is a double covering of the projective plane. The Theorem states that $
\epsilon(\pi^*L,x)=1$, when $x$ is a very general point on the ramification divisor.
Let $2b$ be the degree of the branch divisor. It is well known that $X$ is one of the following surfaces:
\begin{enumerate}

\item If $b=1$, $X\cong \P^1\times \P^1$. The divisor $\pi^*L$ is of type $(1,1)$ on $X$. Its Seshadri constant at any point is $1$.

\item If $b=2$, then $X$ is a Del Pezzo surface. In particular, $X$ is the blowing up of $\P^2$ at $7$ general points. The divisor $\pi^*L$ corresponds to the anticanonical divisor $-K_X$ of $X$. If $x$ is a generic point on $X$ then (see Broustet):
$$
\epsilon(\pi^*L,x)=\frac{4}{3}.
$$
However, when $x$ is on the ramification divisor $\epsilon(\pi^*L,x)=1$.

\item If $b=3$, then $X$ is a $K3$ surface. We prove that $\epsilon(\pi^*L,x)=1$ at a very general point $x$ on the ramification divisor. The same result is a direct consequence of the Example 2.3 of \cite{BaRoSz00}. There, it is proved  that $d\pi^*L$ generates exactly $d$-jets at $x$. The relation between the Seshadri constants and generation of jets provides the desired value of $\epsilon(\pi^*L,x)$ (see Proposition 1.1 of \cite{Ba99}).

\item If $b>3$, then $X$ is a surface of general type. By a result of Buium (see \cite{Bu83}), when the branch divisor is generic, $X$ has Picard number $1$. Thus, we conclude that the global Seshadri constant of any line bundle on a general double cover of the projective plane is rational.

\end{enumerate}

\end{rem}

\section{Relation with the Nagata conjecture.}

Consider the problem of the existence of plane algebraic curves of given degree and with singularities of prescribed order in points in general position (see \cite{StSz04} for a survey on this topic). The Nagata conjecture says:

\begin{conj}[Nagata conjecture]

Let $P_1,\ldots,P_n$ be $n\geq 10$ be general points in $P^2$ and let $k_1,\ldots,k_n$ be fixed non-negative integers. If $C\subset \P^2$ is a curve of degree $d$ such that $mult_{P_i}(C)\geq k_i$ then:
$$
d\geq\frac{1}{\sqrt{n}}\sum_{i=1}^nk_i.
$$
\end{conj}

This can be reformulated in the language of Seshadri constants. If $(Y,L)$ is a polarized variety, $\epsilon(Y;n)$ denotes the multiple Seshadri constant of $L$ at $n$ general points.

\begin{conj}[Nagata conjecture via Seshadri constants]

If $n\geq 9$, the Seshadri constant $\epsilon(\Te_{P^2}(1); n)$ is maximal, that is,
$$
\epsilon(\Te_{P^2}(1); n)=\frac{1}{\sqrt{n}}.
$$
\end{conj}

Let us see the relation between our result and the Nagata conjecture. Let $\pi:X\lrw Y$  be a $n$-cyclic covering of a smooth surface $Y$. Let $L$ be an ample divisor on $Y$. Let $P$ be a generic point on the ramification divisor. The point $Q=\pi(P)\in Y$ is at the branch divisor $B$. Let us suppose that the Seshadri constant $\epsilon(\pi^*L,P)$ does not reach the expected value. In this case, we have seen that the Seshadri constant is given by  a divisor $D=\pi^*C$, with $C\in Div(Y)$. The multiplicity of $D$ at $P$ determines the infinitesimal behavior of $C$ at $Q$. This can be expressed in the language of clusters in the following way (see \cite{Ca90}, \cite{Ca00}).

We consider the following cluster:
$$
\begin{array}{l}
{Y_0=Y,\quad Q_1=Q, \quad B_0=B;}\\
{Y_1=Bl(Y_0,\{Q_1\}),\quad  B_1=B_0',\quad Q_2=E_1\cap B_1;}\\
{Y_2=Bl(Y_1,\{Q_2\}),\quad  B_2=B_1',\quad Q_3=E_2\cap B_2;}\\
{\ldots}\\
{Y_{n-1}=Bl(Y_{n-2},\{Q_{n-1}\}),\quad  B_{n-1}=B_{n-2}',\quad Q_n=E_{n-1}\cap B_{n-1};}\\
\end{array}
$$
where $E_i$ denotes the exceptional divisor of the blowing up of $Y_{i-1}$ and $B_i$ is the strict transform of $B_{i-1}$.

\begin{prop} \label{clusters}

With the previous notation if $D=\pi*C$ is a divisor passing through $P$ with multiplicity $m$ then $C$ is a divisor of $Y$ passing through the cluster $(Q_1,\ldots,Q_n)$ with multiplicities $(m_1,\ldots,m_n)$ and $m=\sum_{i=1}^n m_i$.

\end{prop}
{\bf Proof:} We consider local coordinates $(u,v)$ at $P$ in $X$ such that $v=0$ is the equation 
of the ramification divisor and $u=0$ is the equation of a $\sigma$-invariant divisor. We also take local coordinates $(x,y)$ in $Y$ such that the local expression of the $n$-covering is:
$$
(u,v)\lrw (x,y)=(u,v^n).
$$
and the equation of the branch divisor is $y=0$. The local expression of the automorphism $\sigma$ is:
$$
(u,v)\lrw (u,\theta v)
$$
where $\theta$ is a primitive $n$-root of unity.

Let $f(u,v)=0$ be the local equation of $D$ at $P$. Since $D$ is a $\sigma$-invariant divisor passing through $P$ with multiplicity $m$, the Taylor series of $f$ is:
$$
f(u,v)=\sum_{p+nq\geq m} a_{pq}u^pv^{nq}.
$$
where there is at least a coefficient $a_{pq}\neq 0$ with $p+nq=m$.

The Taylor series of of the curve $C$ such that $D=\pi^*C$ is:
$$
\sum_{p+nq\geq m} a_{pq}x^py^{q}.
$$
Let $m_1=mult_{Q_1}(C)$. It is clear that:
$$
m_1=min\{p+q\}\qquad \hbox{ with } p+nq\geq m\hbox{ and } a_{pq}\neq 0
$$
so
$$
1\leq m_1\leq m.
$$

We make the blowing up $Y_1$ of $Y$ at $Q_1$. We take new coordinates:
$$
x=x_1; \qquad y=x_1y_1.
$$
Now the local equation of the strict transform $B_1$ of the branch divisor $B$ is $y_1=0$. The equation of the exceptional divisor $E_1$ is $x_1=0$. The local expression of the total transform of $C$ is:
$$
\sum_{p+nq\geq m} a_{pq}x_1^{p+q}y_1^{q},
$$
and the expression of the strict transform $C_1$ of $C$ is:
$$
\sum_{p+nq\geq m} a_{pq}x_1^{p+q-m_1}y_1^{q}
$$
If we take $p_1=p+q-m_1$, we have:
$$
\sum_{p_1+(n-1)q\geq m-m_1} a^{(1)}_{p_1q}x_1^{p_1}y_1^{q}
$$
where there is at least a coefficient $a^{(1)}_{p_1q}\neq 0$ with $p_1+(n-1)q=m-m_1\geq 0$.

We can repeat this construction, taking the blowing up of the surfaces $Y_{k-1}$ at the points $Q_k$. In particular, when we construct the surface $Y_{k}=Bl(Y_{k-1},Q_k)$ the local expression of the strict transform $C_k$ of $C_{k-1}$ is:
$$
\sum_{p_{k}+(n-k)q\geq m-m_1-\ldots-m_{k}} a^{(k)}_{p_{k}q}x_1^{p_{k}}y_1^{q}
$$
where there is at least a coefficient $a^{(k)}_{p_{k}q}\neq 0$ with $p_{k}+(n-k)q=m-m_1-\ldots-m_{k}\geq 0$. Note that:
$$
m_{k+1}=mult_{Q_{k+1}}(C_{k})=min\{p_k+q\}
$$
with $p_{k}+(n-k)q\geq m-m_1-\ldots-m_{k}$ and $a^{(k)}_{p_kq}\neq 0$. So
$$
j\leq m_{k+1}\leq m-m_1-\ldots-m_{k}.
$$
where $j=0$ if $m=m_1+\ldots+m_k$ and $j=1$ in other case.

We have two possibilities:

\begin{enumerate}

\item If $m_{n-1}> 0$, then we know that there is at least a coefficient $a^{(n-1)}_{p_{n-1}q}\neq 0$ with $p_{n-1}+q=m-m_1-\ldots-m_{n-1}\geq 0$. Moreover:
$$
m_{n}=mult_{Q_{n}}(C_{n-1})=min\{p_{n-1}+q\}
$$
with $p_{n-1}+q\geq m-m_1-\ldots-m_{n-1}$ and $a^{(n-1)}_{p_{n-1}q}\neq 0$. We deduce that
$$
m_n=m-m_1-\ldots-m_{n-1}.
$$

\item If $m_{n-1}=0$ then $m=m_1+\ldots+m_{n-2}$. The curve $C_{n-1}$ does not pass through the point $Q_{n-1}$ and if we continue with the process we obtain $m=0$. \qed

\end{enumerate}

Now, we can deduce directly the following results:

\begin{teo}

Let $\pi:X\lrw Y$  be a branched $n$-cyclic covering of a smooth surface $Y$. Let $L$ be an ample divisor of $NS(Y)$. Then, if $P_1,\ldots,P_r$ are different points in the branch divisor of $\pi$:
$$
\epsilon(\pi^*L; P_1,\ldots,P_r)\leq n\epsilon(L; nr)
$$
where $\epsilon(L; nr)$ is the Seshadri constant of $L$ at $nr$ very general points. \	qed

\end{teo}

\begin{cor}

Let $\pi:X\lrw Y$  be a branched $n$-cyclic covering of a smooth surface $Y$. Let $L$ be an ample divisor of $NS(Y)$. Then, if $P$ is a point at the branch divisor of $\pi$:
$$
\epsilon(\pi^*L,P)\leq n\epsilon(L; n)
$$
where $\epsilon(L; n)$ is the Seshadri constant of $L$ at $n$ very general points. \qed

\end{cor}

In particular, when the surface $Y$ is the projective plane $\P^2$, we see that the problem of finding a Seshadri exceptional curves for $\pi^*L$ on $X$ is a special case of the general problem of finding exceptional curves on $\P^2$ (see \cite{Ro01}, \cite{Ro04}). If $2\leq n\leq 9$ the exceptional curves described on the Table \ref{tabla1} corresponds to the the exceptional curves providing the  $n$–-tuple Seshadri constants on $\P^2$ (see Example 2.4 of \cite{StSz04}). When $n>9$ the Nagata conjecture is not solved, except if $n$ is a square. This is agree with the Theorem \ref{aproximacion}.

As a consequence of this discussion, we ask for the following Question:

\begin{question}

Let $\pi:X\lrw \P^2$ be a $n$--cyclic covering of the projective plane, with generic branch divisor of degree $nb$ and $n>9$. Is the Seshadri constant of $\pi^*\Te_{P^2}(1)$ at a very general point $\eta$ maximal?.
$$
\epsilon(\pi^*\Te_{P^2}(1),\eta)\stackrel{?}{=}\sqrt{n}.
$$

\end{question}

\begin{rem}

The points $P_1,\ldots,P_n$ of the cluster depends on the branch divisor. Therefore, the condition of ``general points'' on the Nagata conjecture is replaced here for the condition of ``generic branch divisor''.

\end{rem}
\begin{rem}

When $n$ is not a square, an affirmative answer to this question will provide an example of surface with irrational Seshadri constant.

\end{rem}

{\bf E-mail:} lfuentes@udc.es

Luis Fuentes Garc\'{\i}a. 

Departamento de M\'etodos Matem\'aticos y Representaci\'on.

E.T.S. de Ingenieros de Caminos, Canales y Puertos. 

Universidad de A Coru\~na. Campus de Elvi\~na. 15192 A Coru\~na (SPAIN)

\end{document}